\documentclass[12pt]{article}
\usepackage{amsfonts}
 \usepackage{amsthm}
 \usepackage{amsmath}
 \usepackage{amscd}
 \usepackage{amssymb}
 \usepackage{latexsym}
 \usepackage{epsfig}

\usepackage[ascii]{inputenc}
\usepackage[T1]{fontenc}
\usepackage{amsmath,amssymb,amsfonts,textcomp,epsfig}
\usepackage{hyperref}
\usepackage{color}
\usepackage{calc}

\graphicspath{{figs/}}

\def\C{{\mathbb C}}

\def\R{{\mathbb R}}

\newtheorem{remark}{Remark}

\newtheorem{theorem}{Theorem} 
\newtheorem{lemma}{Lemma}
\newtheorem{proposition}{Proposition}
\newtheorem{corollary}{Corollary}

\title{Cubic perturbations of elliptic Hamiltonian vector fields of degree three}
 \author{Lubomir Gavrilov \\
 \normalsize \it Institut de Math\'{e}matiques de Toulouse,
 \normalsize \it Universit\'{e} de Toulouse\\
   \normalsize \it F-31062 Toulouse, \normalsize \it   France
     \\ Iliya D. Iliev\\
\normalsize \it Institute of Mathematics, Bulgarian Academy of Sciences\\
\normalsize \it Bl. 8, 1113 Sofia, Bulgaria  }
\begin{document}
\maketitle
\begin{abstract}
The purpose of the present paper is to  study the limit cycles of  one-parameter perturbed plane Hamiltonian vector field $X_\varepsilon$
$$
X_\varepsilon : \left\{
\begin{array}{llr}
\dot{x}=\;\; H_y+\varepsilon f(x,y)\\
\dot{y}=-H_x+\varepsilon g(x,y),
\end{array} \;\;\;\;\; H~=\frac{1}{2} y^2~+U(x)
 \right.
$$
which bifurcate from the period annuli of $X_0$ for sufficiently small $\varepsilon$.  Here 
 $U$ is a univariate polynomial of degree four without symmetry, and $f, g$  are \emph{arbitrary} cubic polynomials in two variables.

We take a period annulus and parameterize the related
displacement map $d(h,\varepsilon)$ by the Hamiltonian value $h$ and by the 
small parameter $\varepsilon$. Let $M_k(h)$ be the $k$-th coefficient in its expansion
with respect to $\varepsilon$. We establish the general form of $M_k$ and study its zeroes. We 
deduce
 that the period annuli of $X_0$ can produce
  for sufficiently small $\varepsilon$,
 at most 5, 7 or 8 zeroes  in the interior eight-loop case,  the saddle-loop case, and the exterior eight-loop case respectively. In the  interior eight-loop case the bound is exact, while in the saddle-loop case we provide
examples of Hamiltonian fields  which produce
6 small-amplitude limit cycles.  Polynomial perturbations of $X_0$ of higher  degrees are also studied.
\end{abstract}

\newpage
\tableofcontents

\section{Introduction}
We consider cubic systems in the plane which are small perturbations of
Hamiltonian systems with a center. Our goal is to estimate the number of limit
cycles produced by the perturbation. The Hamiltonians we consider have the form
$H=y^2+U(x)$ where $U$ is a polynomial of degree 4. In this paper we exclude
from consideration the four symmetric Hamiltonians $H=y^2+x^2\pm x^4$,
$H=y^2-x^2+x^4$ and $H=y^2+x^4$ because they require a special treatment.
Therefore, one can use the following normal form of the Hamiltonian
\begin{equation}
\label{ham}
H=\frac12y^2+\frac12x^2-\frac23x^3+\frac{a}{4}x^4,\quad a\neq 0,\frac89.
\end{equation}
An easy observation shows that the following four topologically different
cases occur:
$$\begin{array}{ll}
a<0 & \mbox{\rm   saddle-loop,}\\
0<a<1 & \mbox{\rm  eight loop,}\\
a=1 & \mbox{\rm  cuspidal loop,}\\
a>1 & \mbox{\rm  global center.}\\
\end{array}$$
There is one period annulus in the saddle-loop and the global center cases,
two annuli in the cuspidal loop case, and three annuli in the eight loop case.
Take small $\varepsilon>0$ and consider the following one-parameter
perturbation of the Hamiltonian vector field associated to $H$:
\begin{equation}
\label{per}
\begin{array}{l}
\dot{x}=H_y+\varepsilon f(x,y),\\
\dot{y}=-H_x+\varepsilon g(x,y),\end{array}
\end{equation}
where $f$ and $g$ are arbitrary cubic polynomials with
coefficients $a_{ij}$ and $b_{ij}$ at $x^iy^j$, respectively.
As well known, if we parameterize the displacement map by the
Hamiltonian level $h$, then the following expansion formula holds
\begin{equation}
\label{exp}
d(h,\varepsilon)=\varepsilon M_1(h)+\varepsilon^2 M_2(h)+\varepsilon^3 M_3(h)
+\ldots, \quad h\in \Sigma
\end{equation}
where $\Sigma$ is an open interval depending on the case and the period annulus we consider.
There is a lot of papers investigating system (\ref{per}), but most of them
deal with $M_1(h)$ only or consider perturbations like $f(x,y)=0$, $g(x,y)=(\alpha_0+\alpha_1x+\alpha_2x^2)y$.
See e.g. the book by Colin Christopher and Chengzhi Li \cite{chli07} for more comments and references. 
In what follows we consider for a first time the full 20-parameter cubic deformation (\ref{per}) of the Hamiltonian system
associated to $H$. We suppose, however, that the arbitrary cubic polynomials $f, g$ do not depend on the small parameter $\varepsilon$. To study the full neighborhood of the Hamiltonian system associated to $H$, 
it is also necessary to allow that $f,g$ depend  analytically on $\varepsilon$. 

Our first goal will be to calculate explicitly the first several coefficients
$M_1$, $M_2$, etc. in (\ref{exp}) and then determine the least integer
$m$ such that system (\ref{per}) becomes integrable provided that the first
$m$ coefficients in (\ref{exp}) do vanish.

Let us rewrite system (\ref{per}) in a Pfaffian form
\begin{equation}
\label{pfaff}
dH=\varepsilon\omega,\quad \omega=g(x,y)dx-f(x,y)dy.
\end{equation}
We first establish that if $M_1(h)\equiv 0$, then one can express the
cubic one-form $\omega$ in the perturbation as
\begin{equation}
\label{one}
\textstyle
\omega=d[Q(x,y)-(\frac{a}{5}\lambda-\frac25\mu)x^5-\frac{a}{6}\mu x^6]
+(\lambda x+\mu x^2)dH
\end{equation}
where $Q(x,y)=\sum_{1\leq i+j\leq 4}q_{ij}x^iy^j$ and $\lambda$, $\mu$ are
parameters. Obviously, there are simple explicit linear formulas connecting
$q_{ij}$, $\lambda$ and $\mu$ to the coefficients of $f$ and $g$. We shall
consider $q_{ij}$, $\lambda$ and $\mu$ as the parameters of the perturbation.

\begin{theorem}\label{t1}
The perturbation $(\ref{per})$, $(\ref{pfaff})-(\ref{one})$ is integrable if
and only if either of the two conditions holds:

\vspace{1ex}
{\em 1)}  $\lambda=\mu=0$;

{\em 2)}  $q_{01}=q_{11}=q_{21}=q_{31}=q_{03}=q_{13}=0$.

\vspace{1ex}
\noindent
In the first case system $(\ref{per})$, $(\ref{pfaff})$ becomes
Hamiltonian and in the second one it becomes time-reversible.  

 If
$M_1(h)=M_2(h)=M_3(h)=M_4(h)\equiv 0$, then the perturbation is integrable.
\end{theorem}

\noindent
When the perturbation is integrable, all coefficients $M_k(h)$ do vanish
in the respective period annulus and the Poincar\'e map is the identity.
When the perturbation is not integrable (that is neither of the conditions
in Theorem \ref{t1} holds), one can prove the following result.
Take an oval $\delta(h)$ contained in the level set $H=h$, $h\in\Sigma$
and define the integrals $$I_k(h)=\oint_{\delta(h)}x^kydx, k=0,1,2,\ldots$$

\begin{theorem}\label{t2}
The first four coefficients $M_k(h)$, $1\leq k\leq 4$ have the form
$$M_k(h)=\alpha_k(h)I_0(h)+\beta_k(h)I_1(h)+\gamma_k(h)I_2(h)$$
where $\alpha_k(h)$, $\beta_k(h)$, $\gamma_k(h)$ are polynomials of degree
at most one. The second coefficient $M_2(h)$ has the maximum possible number
of zeroes in $\Sigma$ among $M_k(h)$.
\end{theorem}

We use the above results in deriving upper bounds for the number of limit cycles 
bifurcating from the open period annuli in the cases when the Hamiltonian has
three real and different critical values. For this, we take a perturbation with 
$M_1(h)\equiv 0$ and $M_2(h)\not\equiv 0$, with all six coefficients 
independently free.  

\begin{theorem}\label{t3} {\em (i)} In the interior eight-loop case, at most five 
limit cycles bifurcate from each one of the annuli inside the loop.

{\em (ii)} In the exterior eight-loop case, at most eight limit cycles bifurcate 
from the annulus outside the loop. 

{\em (iii)} In the saddle-loop case, at most seven limit cycles bifurcate from the unique 
period annulus.   
\end{theorem}

\noindent
The proof is based on a refinement of Petrov's method which we apply to 
the much more general case when the coefficients in $M_k(h)$ are polynomials 
of arbitrary degree $n$, thus $M_k(h)$ being an element of a module of 
dimension $3n+3$.

\begin{theorem}\label{t4} Let the coefficients $\alpha_k(h)$, $\beta_k(h)$ and $\gamma_k(h)$
in the expression of $M_k(h)$ be polynomials of degree $n$ with real coefficients.
Then $M_k(h)$ has in the respective interval $\Sigma$ at most $3n+2$ zeroes in the interior 
eight-lop case, at most $4n+4$ in the exterior eight-loop case, and at most $4n+3$ zeroes in
the saddle-loop case.  
\end{theorem}

\noindent
In order to demonstrate that Chebyshev's property (no more zeroes than the dimension minus one) would
not also hold in the saddle-loop case,  we provide an estimate from below for the number of bifurcating 
small-amplitude limit cycles around the center at the origin which concerns all Hamiltonian parameters $a\neq 0, \frac89$. 

\begin{theorem}\label{t5} For $a$ close to $-\frac83$, function $M_1(h)$ can produce 
four small limit cycles around the origin. For $a$ close to $-\frac89$, function $M_2(h)$ 
can produce six such limit cycles. For all other values of $a\in\R$, the number of small
limit cycles produced by the function $M_k(h)$ equals its dimension minus one. 
\end{theorem}

\noindent
The limit cycle in addition in the saddle-loop case is obtained by moving slightly the Hamiltonian
parameter $a$ in appropriate direction from the respective fraction.  

The paper is organized as follows. At the beginning, we compute explicitly the coefficients $M_k$ for $k=1, 2, 3, 4$.
 It is easily seen that for each $k$ they form a set which is
\begin{itemize}
\item a vector space of dimension four, for $k=1$
\item a vector space of dimension six, for $k=2$
\item a union of three distinct five-dimensional vector spaces, for $k=3$
\item a union of three distinct straight lines, for $k=4$,
\end{itemize}
and when $M_1=M_2=M_3=M_4=0$, then the perturbation becomes integrable. The function $M_2$ takes therefore the form
\begin{equation}
\label{m22}
M_2(h) = \alpha_1I_0(h)+\beta_1I_1(h)+\gamma_1I_2(h)
\end{equation}
where $\alpha, \beta, \gamma$ are arbitrary linear functions in $h$. 

Next, considering the generalized
situation when $M_2$ is a function of the form (\ref{m22}) in which $\alpha, \beta, \gamma$ are arbitrary degree $n$ polynomials in $h$, we establish that $M_2$ would have at most
$3n+2$ zeroes in the interior eight-loop case,
$4n+4$ zeroes in the exterior eight-loop case,
$4n+3$ zeroes in the saddle-loop case.
We apply these results to our problem by taking $n=1$. Finally, we provide
examples of Hamiltonian fields in the saddle-loop case which produce
4 and 6 small-amplitude limit cycles, respectively when $M_1\not\equiv 0$,
and $M_1\equiv 0$ but $M_2\not\equiv 0$. For all other cases, the number of such small-amplitude 
limit cycles is less than the respective dimension.

\section{Calculation of the coefficients $M_k(h)$}

In this section we are going to calculate the first four coefficients in
(\ref{exp}). We use the recursive procedure proposed by Fran\c{c}oise
\cite{fran96}, see also \cite{ilie96}, \cite{ilie98a}.

\subsection{The coefficient $M_1(h)$}
We begin with the easy calculation of $M_1(h)$.

\begin{proposition}\label{p1}
{\rm (i)} The function $M_1(h)$ has the form
\begin{equation}\label{m1}
M_1(h)=\alpha_1I_0(h)+\beta_1I_1(h)+\gamma_1I_2(h),
\end{equation}
where $\alpha_1$ is a first-degree polynomial in $h$ and $\beta_1$, $\gamma_1$
are constants, depending on the perturbation.

\vspace{1ex}
\noindent
{\rm (ii)} If $M_1(h)\equiv 0$, then one can rewrite the one-form $\omega$ as
$(\ref{one})$ where $Q$ is a polynomial of degree four without constant term
and $\lambda$, $\mu$ are constant parameters.
\end{proposition}

\vspace{1ex}
\noindent
{\bf Proof.} By a simple calculation, one can rewrite $\omega$ in the form
$\omega=dQ(x,y)+yq(x,y)dx$ with $Q$ and $q$ certain polynomials of degree
4 and 2, respectively. Denote for a moment by $c_{ij}$ the coefficient in
$q$ at $x^iy^j$. Then
$$yq(x,y)dx=(c_{01}+c_{11}x)y^2dx+(c_{00}+c_{10}x+c_{20}x^2)ydx+c_{02}y^3dx.$$
Next,
$y^3dx=(2H-x^2+\frac43x^3-\frac{a}{2}x^4)ydx=(2H-x^2)ydx
+yd(\frac13x^4-\frac{a}{10}x^5).$
Using the identity
$
\frac13x^4-\frac{a}{10}x^5=\frac{4}{15a}H-\frac25xH+(\frac{8}{45a}+\frac15)x^3
-\frac{2}{15a}x^2-\frac{2}{15a}y^2+\frac15xy^2$
we derive the equation
\begin{equation}\label{rec}
\textstyle
y^3dx= d(\frac17xy^3-\frac{2}{21a}y^3)+(\frac{2}{7a}-\frac37x)ydH
+[\frac{12}{7}H-\frac{2}{7a}x+(\frac{4}{7a}-\frac37)x^2]ydx.
\end{equation}
Replacing in the formula above and taking into account that
$M_1(h)=\oint_{\delta(h)}\omega=\oint_{\delta(h)}yq(x,y)dx$, one obtains
formula (\ref{m1}) with
$$\textstyle \alpha_1=c_{00}+\frac{12}{7}c_{02}h,\quad
\beta_1=c_{10}-\frac{2}{7a}c_{02},\quad
\gamma_1=c_{20}+(\frac{4}{7a}-\frac37)c_{02}.$$
Now, $M_1(h)\equiv 0$ is equivalent to $c_{00}=c_{10}=c_{20}=c_{02}=0$
(see Corollary \ref{star} below)
and $\omega$ becomes  $\omega=dQ-\frac12y^2(\lambda+2\mu x)dx$
where $\lambda=-2c_{01}$, $\mu=-c_{11}$. On the other hand (modulo terms $dQ$)
$$\begin{array}{l}
-\frac12y^2(\lambda+2\mu x)dx=(\lambda x+\mu x^2)d(H-\frac12x^2+\frac23x^3
-\frac{a}{4}x^4)\\
=(\lambda x+\mu x^2)dH+(\lambda x+\mu x^2)(-x+2x^2-ax^3)dx\\
= (\lambda x+\mu x^2)dH+d(-\frac{a}{5}\lambda x^5+\frac25\mu x^5
-\frac{a}{6}\mu x^6). \end{array}$$
Proposition \ref{p1} is proved.  $\Box$

\vspace{1ex}
\noindent

\subsection{The coefficient $M_2(h)$}
By (\ref{one}), if $\lambda=\mu=0$, then the perturbation is Hamiltonian
and all coefficients $M_k$ do vanish. We will assume below that $\lambda$
and $\mu$ are not both zero. Then the calculation of $M_2(h)$ makes sense.
Denote by $q_{ij}$ the coefficient at $x^iy^j$ in $Q$. Below, we split $Q$
into an odd and even part $Q=Q_1+Q_2$ with respect to $y$.

\begin{proposition}\label{p2}
{\rm (i)} If $M_1(h)\equiv 0$, then the function $M_2(h)$ has the form
\begin{equation}\label{m2}
M_2(h)=\alpha_2I_0(h)+\beta_2I_1(h)+\gamma_2I_2(h),
\end{equation}
where $\alpha_2$, $\beta_2$ and $\gamma_2$
are first-degree polynomials in $h$ with coefficients depending on the
perturbation.

\vspace{1ex}
\noindent
{\rm (ii)} If $M_1(h)=M_2(h)\equiv 0$, then the odd part of $Q(x,y)$ becomes:

\vspace{1ex}
{\rm (a)}  $Q_1=q_{11}(x-2x^2+ax^3)y,\quad$ if $\mu=0$;

\vspace{1ex}
{\rm (b)}  $Q_1=-\frac12q_{11}(1-2x+ax^2)y,\quad$ if $\lambda=0$;

\vspace{1ex}
{\rm (c)}  $Q_1=q_{11}(x+\frac{a\lambda}{2\mu}x^2)y,\quad$ if $a\leq 1$ and
$a\lambda^2+4\lambda\mu+4\mu^2= 0$;

\vspace{1ex}
{\rm (d)}  $Q_1=0,\quad$ if $\lambda\mu\neq 0$ and
$a\lambda^2+4\lambda\mu+4\mu^2\neq 0$.
\end{proposition}

\vspace{1ex}
\noindent
{\bf Proof.} As well known, the second coefficient in (\ref{exp}) is obtained
by integrating the one-form $\omega_2=(\lambda x+\mu x^2)\omega$, that is
$$M_2(h)=\oint_{\delta(h)}\omega_2=\oint_{\delta(h)}(\lambda x+\mu x^2)dQ(x,y)=
-\oint_{\delta(h)} (\lambda +2\mu x) Q_1(x,y)dx$$
$$=-\oint_{\delta(h)} (\lambda +2\mu x)[(q_{01}+q_{11}x+q_{21}x^2+q_{31}x^3)y
+(q_{03}+q_{13}x)y^3]dx.$$
Next, multiplying (\ref{rec}) by $x$ and expressing the first term on the
right-hand side in a proper form, we obtain identity
\begin{equation}\label{rec1}
\begin{array}{l}
xy^3dx=d(\frac18x^2y^3-\frac{1}{14a}xy^3-\frac{1}{126a^2}y^3)
+(\frac{1}{42a^2}+\frac{3}{14a}x-\frac38x^2)ydH\\
+[(\frac{1}{7a}+\frac32x)H-\frac{1}{42a^2}x+(\frac{1}{21a^2}-\frac{2}{7a})x^2
+(\frac{1}{2a}-\frac38)x^3]ydx.
\end{array}
\end{equation}
In a similar way, multiplying (\ref{rec1}) by $x$, we get
\begin{equation}\label{rec2}
\begin{array}{l}
x^2y^3dx=d(\frac19x^3y^3-\frac{1}{18a}x^2y^3-\frac{2}{189a^2}xy^3
-\frac{2}{1701a^3}y^3)\\
+(\frac{2}{567a^3}+\frac{2}{63a^2}x+\frac{1}{6a}x^2-\frac13x^3)ydH
+[(\frac{4}{189a^2}+\frac{2}{9a}x+\frac43x^2)H\\
-\frac{2}{567a^3}x+(\frac{4}{567a^3}-\frac{8}{189a^2})x^2
+(\frac{2}{27a^2}-\frac{5}{18a})x^3+(\frac{4}{9a}-\frac13)x^4]ydx.
\end{array}
\end{equation}
Replacing the values from (\ref{rec}), (\ref{rec1}) and (\ref{rec2})
in the above formula of $M_2(h)$, we obtain
$$
M_2(h)=-[q_0I_0(h)+q_1I_1(h)+q_2I_2(h)+q_3I_3(h)+q_4I_4(h)]
$$
where
$$\begin{array}{rl}
q_0=& \lambda q_{01}+[\frac{12}{7}\lambda q_{03}
+\frac{1}{7a}(\lambda q_{13}+2\mu q_{03})+\frac{8}{189a^2}\mu q_{13}]h,\\
q_1=& \lambda q_{11}+2\mu q_{01}-\frac{2}{7a}\lambda q_{03}
-\frac{1}{42a^2}(\lambda q_{13}+2\mu q_{03})-\frac{4}{567a^3}\mu q_{13}\\
    &+[\frac32\lambda q_{13}+3\mu q_{03}+\frac{4}{9a}\mu q_{13}]h,\\
q_2=& \lambda q_{21}+2\mu q_{11} +(\frac{4}{7a}-\frac37)\lambda q_{03}
+(\frac{1}{21a^2}-\frac{2}{7a})(\lambda q_{13}+2\mu q_{03})\\
 & +(\frac{8}{567a^3}-\frac{16}{189a^2})\mu q_{13}+\frac83\mu q_{13}h,\\
q_3=&\lambda q_{31}+2\mu q_{21}+(\frac{1}{2a}-\frac38)
(\lambda q_{13}+2\mu q_{03})+(\frac{4}{27a^2}-\frac{5}{9a})\mu q_{13},\\
q_4= &2\mu q_{31}+(\frac{8}{9a}-\frac23)\mu q_{13}.
\end{array}$$
In order to remove integrals $I_3, I_4$, we use the identity
$$\oint_{\delta(h)}(x^k U'+{\textstyle\frac23}kx^{k-1}U)ydx=0,\quad
\textstyle {U=h-\frac12x^2+\frac23x^3-\frac14ax^4}$$
which is equivalent to
\begin{equation}\label{reck}
\textstyle\frac{k+6}{6}aI_{k+3}=\frac{4k+18}{9}I_{k+2}-\frac{k+3}{3}I_{k+1}
+\frac{2k}{3}hI_{k-1}.
\end{equation}
Used with $k=0,1,2$, this relation yields
\begin{equation}\label{I345}
\begin{array}{l}
I_3=\frac{2}{a}I_2-\frac{1}{a}I_1,\\
I_4=(\frac{88}{21a^2}-\frac{8}{7a})I_2-\frac{44}{21a^2}I_1+\frac{4}{7a}hI_0,\\
I_5=(\frac{572}{63a^3}-\frac{209}{42a^2})I_2-
(\frac{286}{63a^3}-\frac{5}{4a^2}-\frac{1}{a}h)I_1+\frac{26}{21a^2}hI_0.
\end{array}
\end{equation}
Replacing, we finally derive formula (\ref{m2}) with
$$\begin{array}{l}
\alpha_2=-q_0-\frac{4}{7a}q_4h,\\[1ex]
\beta_2=-q_1+\frac{1}{a}q_3 +\frac{44}{21a^2}q_4,\\[1ex]
\gamma_2=-q_2-\frac{2}{a}q_3+(\frac{8}{7a}-\frac{88}{21a^2})q_4.
\end{array}$$
Then $M_2(h)\equiv 0$ is equivalent to $\alpha_2=\beta_2=\gamma_2=0$
(see Corollary \ref{star} below). Taking the coefficients at $h$ zero, we obtain that
either $\mu=q_{03}=q_{13}=0$ or $\mu\neq 0$ and $q_{31}=q_{03}=q_{13}=0$.
In the first case, $\lambda\neq 0$ and taking the coefficients at 1 zero, we
easily obtain $q_{01}=0$, $q_{21}=-2q_{11}$, $q_{31}=aq_{11}$ which is case (a).
In the second case above, if $\lambda=0$, then one easily obtains
$q_{01}=-\frac12q_{11}$, $q_{21}=-\frac{a}{2}q_{11}$ which is case (b).
If $\lambda\neq 0$, then taking the coefficients at 1 zero yields
$q_{01}=0$ and equations $-\lambda q_{11}+\frac{2}{a}\mu q_{21}=0$,
$-2\mu q_{11}-(\lambda+\frac{4}{a}\mu) q_{21}=0$. Provided that
$a\lambda^2+4\lambda\mu+4\mu^2=0$ (it is possible for $a\leq 1$ only),
one has $q_{21}=\frac{a\lambda}{2\mu}q_{11}$ which is case (c).
Otherwise, one obtains $q_{11}=q_{21}=0$ which is case (d). Proposition
\ref{p2} is proved.   $\Box$

\vspace{1ex}
\noindent
\subsection{The coefficient $M_3(h)$}
It turns out that if $Q_1=0$ then the perturbation is integrable. This is
because the perturbed system (\ref{per}) becomes time-reversible in this case.
Below we are going to consider the three cases (a), (b), (c) when
$q_{11}\neq 0$. For them, the next coefficient $M_3(h)$ in (\ref{exp}) should
be calculated. For this purpose, we need to express the one-form
$\omega_2=(\lambda x+\mu x^2)\omega$ as $dS_2+s_2dH$ and then integrate
the one-form $\omega_3=s_2\omega$.

\begin{proposition}\label{p3} Assume that $q_{11}\neq 0$.

\noindent
{\rm (i)} If $M_1(h)=M_2(h)\equiv 0$, then the function $M_3(h)$ has the form
\begin{equation}\label{m3}
M_3(h)=\alpha_3I_0(h)+\beta_3I_1(h)+\gamma_3I_2(h),
\end{equation}
where $\alpha_3$, $\beta_3$ are first-degree polynomials in $h$ with
coefficients depending on the perturbation and $\gamma_3$ is a constant.

\vspace{1ex}
\noindent
{\rm (ii)} If $M_1(h)=M_2(h)=M_3(h)\equiv 0$, then the even part of $Q(x,y)$
becomes

\vspace{1ex}
$\!\! Q_2=q_{20}x^2-(\frac43q_{20}+\frac13\lambda)x^3+(\frac{a}{2}q_{20}
+\frac12\lambda-\frac14\mu)x^4+(q_{02}-\frac13\lambda x-\frac13\mu x^2)y^2
+q_{04}y^4$

\vspace{1ex}
\noindent
where $\mu=0$ in case {\em (a)}, $\lambda=0$ in case {\em (b)} and
$a\lambda^2+4\lambda\mu+4\mu^2=0$ in case {\em (c)}.
\end{proposition}

\vspace{1ex}
\noindent
{\bf Proof.} To find $s_2$, is suffices to perform the calculations modulo exact
forms. Let us handle first case (a). By (\ref{one}) one obtains (neglecting the
exact forms)
$$\omega_2=\lambda xd[(q_{02}+q_{12}x+q_{22}x^2)y^2+q_{04}y^4+
q_{11}(x-2x^2+ax^3)y]+\lambda^2x^2dH.$$
Then $xdq_{11}(x-2x^2+ax^3)y=-q_{11}(x-2x^2+ax^3)ydx=
-q_{11}yd(H-\frac12y^2)=-q_{11}ydH$.
Similarly,
$$\begin{array}{l}
xd(q_{02}+q_{12}x+q_{22}x^2)y^2=
2xd(q_{02}+q_{12}x+q_{22}x^2)H\\
=-2(q_{02}+q_{12}x+q_{22}x^2)Hdx
=2(q_{02}x+\frac12q_{12}x^2+\frac13q_{22}x^3)dH.\end{array}$$
Finally,
$$\begin{array}{l}
xdq_{04}y^4=
q_{04}xd[4H^2-4H(x^2-\frac43x^3+\frac{a}{2}x^4)] =8q_{04}xHdH\\
+4q_{04}H(x^2-\frac43x^3+\frac{a}{2}x^4)dx
=4q_{04}(2xH-\frac13x^3+\frac13x^4-\frac{a}{10}x^5)dH.
\end{array}$$
Summing up all terms together, we obtain for case (a)
$$\textstyle
s_2=\lambda^2x^2-\lambda q_{11}y
+2\lambda(q_{02}x+\frac12q_{12}x^2+\frac13q_{22}x^3)
+4\lambda q_{04}(2xH-\frac13x^3+\frac13x^4-\frac{a}{10}x^5).
$$

In a similar way, we consider (b). In this case,
$$\textstyle
\omega_2=\mu x^2d[(q_{02}+q_{12}x+q_{22}x^2)y^2+q_{04}y^4-\frac12
q_{11}(1-2x+ax^2)y]+\mu^2x^4dH.$$
Then $-\frac12x^2d(1-2x+ax^2)y=(x-2x^2+ax^3)ydx=ydH$,
$$\begin{array}{l}
x^2d(q_{02}+q_{12}x+q_{22}x^2)y^2=
2x^2d(q_{02}+q_{12}x+q_{22}x^2)H\\
=-4(q_{02}x+q_{12}x^2+q_{22}x^3)Hdx
=4(\frac12q_{02}x^2+\frac13q_{12}x^3+\frac14q_{22}x^4)dH,\end{array}$$
$$\begin{array}{l}
x^2dy^4=
x^2d[4H^2-4H(x^2-\frac43x^3+\frac{a}{2}x^4)] =8x^2HdH\\
+8H(x^3-\frac43x^4+\frac{a}{2}x^5)dx
=8(x^2H-\frac14x^4+\frac{4}{15}x^5-\frac{a}{12}x^6)dH.
\end{array}$$
Summing up all needed terms, we obtain in case (b) the formula
$$\textstyle
s_2=\mu^2x^4+\mu q_{11}y
+4\mu(\frac12q_{02}x^2+\frac13q_{12}x^3+\frac14q_{22}x^4)
+8\mu q_{04}(x^2H-\frac14x^4+\frac{4}{15}x^5-\frac{a}{12}x^6).
$$
Finally, in case (c) we have $a\lambda^2+4\lambda\mu+4\mu^2=0$ and
$$\textstyle
\omega_2=(\lambda x+\mu x^2)d[(q_{02}+q_{12}x+q_{22}x^2)y^2+q_{04}y^4+
q_{11}(x+\frac{a\lambda}{2\mu}x^2)y]+(\lambda x+\mu x^2)^2dH.$$
As above,
$$\begin{array}{l}
(\lambda x+\mu x^2)dq_{11}(x+\frac{a\lambda}{2\mu}x^2)y=
-q_{11}(x+\frac{a\lambda}{2\mu}x^2)(\lambda+2\mu x)ydx\\
=-q_{11}yd(\frac{\lambda}{2}x^2+\frac{a\lambda^2+4\mu^2}{6\mu}x^3+
\frac{a\lambda}{4}x^4)=-\lambda q_{11}yd(H-\frac12y^2)=-\lambda q_{11}ydH,
\end{array}$$
$$\begin{array}{l}
(\lambda x+\mu x^2)d(q_{02}+q_{12}x+q_{22}x^2)y^2
=-2H(q_{02}+q_{12}x+q_{22}x^2)(\lambda+2\mu x)dx\\
=[2\lambda q_{02}x+(\lambda q_{12}+2\mu q_{02})x^2
+\frac23(\lambda q_{22}+2\mu q_{12})x^3+\mu q_{22}x^4]dH,\end{array}$$
$$\begin{array}{l}
(\lambda x+\mu x^2)dq_{04}y^4=
q_{04}(\lambda x+\mu x^2)d[4H^2-4H(x^2-\frac43x^3+\frac{a}{2}x^4)]\\
=8q_{04}(\lambda x+\mu x^2)HdH+4q_{04}H(x^2-\frac43x^3+\frac{a}{2}x^4)(\lambda+2\mu x)dx\\
=4q_{04}[2(\lambda x+\mu x^2)H-\frac13\lambda x^3+(\frac13\lambda-\frac12\mu)x^4
-(\frac{a}{10}-\frac{8}{15}\mu)x^5-\frac{a}{6}\mu x^6]dH.
\end{array}$$
Summing up all terms, we obtain in case (c) the respective formula
$$\begin{array}{l}
s_2=(\lambda x+\mu x^2)^2-\lambda q_{11}y+
2\lambda q_{02}x+(\lambda q_{12}+2\mu q_{02})x^2
+\frac23(\lambda q_{22}+2\mu q_{12})x^3\\
+\mu q_{22}x^4+4q_{04}[2(\lambda x+\mu x^2)H-\frac13\lambda x^3+
(\frac13\lambda-\frac12\mu)x^4-(\frac{a}{10}\lambda-\frac{8}{15}\mu)x^5
-\frac{a}{6}\mu x^6].\end{array}$$

In order to calculate $M_3$ at once for all three cases (a), (b), (c), we
shall use the formula of $s_2$ for case (c) from which the other two cases
are obtained by taking $\mu$ or $\lambda$ zero. Indeed, let us denote by
$s_2^0$ the even part of $s_2$ with respect to $y$. Then
$s_2=\kappa y+s_2^0$ where $\kappa=-\lambda q_{11}$ in cases (a), (c)
and $\kappa=\mu q_{11}$ in case (b). Then
$$M_3(h)=\oint_{\delta(h)}s_2\omega=\oint_{\delta(h)}\kappa yd[Q_2
+({\textstyle\frac25}\mu-{\textstyle\frac{a}{5}}\lambda)x^5\
-{\textstyle\frac{a}{6}}\mu x^6]+\oint_{\delta(h)}s_2^0dQ_1=I+J.$$
We further have
$$I=\kappa\oint_{\delta(h)}[(q_{10}+2q_{20}x+3q_{30}x^2+4q_{40}x^3+
(2\mu-a\lambda)x^4-a\mu x^5)y+(\textstyle\frac13q_{12}
+\textstyle\frac23 q_{22}x)y^3]dx$$
$$=\kappa(q_0I_0+q_1I_1+q_2I_2+q_3I_3+q_4I_4+q_5I_5)$$
with
$$\begin{array}{l}
q_0=q_{10}+(\frac47q_{12}+\frac{2}{21a}q_{22})h,\\
q_1=2q_{20}-\frac{2}{21a}q_{12}-\frac{1}{63a^2}q_{22}+q_{22}h,\\
q_2=3q_{30}+(\frac{4}{21a}-\frac17)q_{12}+(\frac{2}{63a^2}-\frac{4}{21a})q_{22},\\
q_3=4q_{40}+(\frac{1}{3a}-\frac14)q_{22},\\
q_4=2\mu-a\lambda,\\
q_5=-a\mu,\end{array}$$
(we used (\ref{rec}) and (\ref{rec1}) as well).
On the other side, integrating by parts one can rewrite $J$ as
$J=-\oint_{\delta(h)}(s_2^0)'Q_1dx =J_1+J_2$
where $J_2$ is the part corresponding to the expression in $s_2^0$
which contains $q_{04}$. Let us first verify that $J_2=0$.
Indeed, one can establish by easy calculations that
$$\begin{array}{l}
4q_{04}[(2\lambda+4\mu x)H-\lambda x^2+
(\frac43\lambda-2\mu)x^3-(\frac{a}{2}\lambda-\frac{8}{3}\mu)x^4
-a\mu x^5]\\
=8q_{04}(\lambda+2\mu x)(H-\frac12x^2+
\frac23x^3-\frac{a}{4}x^4)
=4q_{04}(\lambda+2\mu x)y^2,\\
-Q_1(\lambda+2\mu x)=\kappa y(x-2x^2+ax^3).\end{array}$$
Hence,
$$J_2=4\kappa q_{04}\oint_{\delta(h)}(x-2x^2+ax^3)y^3dx=
4\kappa q_{04}\oint_{\delta(h)}y^3d(H-{\textstyle \frac12}y^2)=0.$$
What $J_1$ concerns, another easy calculation shows that
$$\begin{array}{l}
-2[(\lambda x+\mu x^2)(\lambda+2\mu x)+\lambda q_{02}+(\lambda q_{12}
+2\mu q_{02})x +(\lambda q_{22}+2\mu q_{12})x^2+2\mu q_{22}x^3]Q_1\\
=2\kappa (x-2x^2+ax^3)[q_{02}+(q_{12}+\lambda)x+(q_{22}+\mu)x^2]y
\end{array}$$
for all three cases. Therefore, by integrating, one obtains
$$J=J_1=\kappa(r_1I_1+r_2I_2+r_3I_3+r_4I_4+r_5I_5)$$
where
$$\begin{array}{l}
r_1=2q_{02},\\
r_2=2\lambda-4q_{02}+2q_{12},\\
r_3=2\mu-4\lambda+2aq_{02}-4q_{12}+2q_{22},\\
r_4=2a\lambda-4\mu+2aq_{12}-4q_{22},\\
r_5=2a\mu+2aq_{22}.
\end{array}$$
Combining with the formula of $I$ and using (\ref{I345}), one obtains
expression (\ref{m3}) with coefficients
$$\begin{array}{l}
\alpha_3=\kappa[q_0+\frac{4}{7a}h(q_4+r_4)+\frac{26}{21a^2}h(q_5+r_5)],\\
\beta_3=\kappa[q_1+r_1-\frac{1}{a}(q_3+r_3)-\frac{44}{21a^2}(q_4+r_4)-
(\frac{286}{63a^3}-\frac{5}{4a^2}-\frac{1}{a}h)(q_5+r_5)],\\
\gamma_3=\kappa[q_2+r_2+\frac{2}{a}(q_3+r_3)+(\frac{88}{21a^2}-\frac{8}{7a})(q_4+r_4)
+(\frac{572}{63a^3}-\frac{209}{42a^2})(q_5+r_5)].
\end{array}$$
It is seen that $\alpha_3$ and $\beta_3$ are first-degree polynomials while
$\gamma_3$ is a constant polynomial. This proves part (i) of the statement.
To prove part (ii), assume that $M_3(h)$ vanishes, which is equivalent to
$\alpha_3=\beta_3=\gamma_3=0$ (see Corollary \ref{star} below). Then by
straightforward calculations one obtains that this is equivalent to
$$\textstyle
q_{10}=0,\;\;\;
q_{30}=-\frac43q_{20}-\frac13\lambda,\;\;\;
q_{40}=\frac{a}{2}q_{20}+\frac12\lambda-\frac14\mu,\;\;\;
q_{12}=-\frac13\lambda,\;\;\;
q_{22}=-\frac13\mu$$
which yields the needed formula of $Q_2$. Proposition \ref{p3} is proved.   $\Box$

\vspace{2ex}
\noindent

\subsection{The coefficient $M_4(h)$}
Replacing the values of the coefficients we just calculated, we obtain
$$\begin{array}{rl}
\omega=&(2q_{20}-\lambda x-\mu x^2)(x-2x^2+ax^3)dx\\
&+d[Q_1+(q_{02}-\frac13\lambda x-\frac13 \mu x^2)y^2+q_{04}y^4]+
(\lambda x+\mu x^2)dH,\\
s_2=&\frac23(\lambda x+\mu x^2)^2+\kappa y+2q_{02}(\lambda x+\mu x^2)\\
&+4q_{04}[2(\lambda x+\mu x^2)H-\lambda(\frac13x^3-\frac13x^4+\frac{a}{10}x^5)
-\mu(\frac12x^4-\frac{8}{15}x^5+\frac{a}{6}x^6)].\end{array}$$

\begin{proposition}\label{p4} Assume that $q_{11}\neq 0$ and
$M_1(h)=M_2(h)=M_3(h)\equiv 0$. Then the function $M_4(h)$ has the form
$$\begin{array}{l}
M_4(h)=\lambda q_{11}^3[2hI_0(h)-(3ah+\frac34-\frac{2}{3a})I_1(h)
+(\frac32-\frac{4}{3a})I_2(h)],\quad \mu=0,\\[1ex]
M_4(h)=-\frac12\mu q_{11}^3[I_0(h)-2I_1(h)+aI_2(h)],\quad \lambda=0,\\[1ex]
M_4(h)=-(\frac{\lambda^2}{\mu^2}+\frac{3\lambda}{2\mu}) q_{11}^3[2\mu I_1(h)
+a\lambda I_2(h)],\quad \lambda\mu\neq 0,\quad
a\lambda^2+4\lambda\mu+4\mu^2= 0.
\end{array}$$
Moreover, $M_4(h)\not\equiv 0$.
\end{proposition}

\vspace{1ex}
\noindent
{\bf Proof.} In what follows, it is useful to introduce notations
$$\begin{array}{l}
A=\lambda(\frac13x^3-\frac13x^4+\frac{a}{10}x^5)+
\mu(\frac12x^4-\frac{8}{15}x^5+\frac{a}{6}x^6),\\
B=\frac12x^2-\frac23x^3+\frac{a}{4}x^4,\qquad L=\lambda x+\mu x^2.
\end{array} $$
Then $dA=2BdL$, $(2q_{20}-L)dB=d[(2q_{20}-L)B+\frac12A]$ and one can rewrite
the expressions of $\omega$ and $s_2$
as follows:
$$\begin{array}{rl}
\omega &=(2q_{20}-L)B'dx+d[Q_1+(q_{02}-\frac13L)y^2+q_{04}y^4]+ LdH,\\
 &=d[Q_1+ (2q_{20}-L)B+\frac12A+(q_{02}-\frac13L)y^2+q_{04}y^4]+ LdH,\\
s_2&=\frac23L^2+\kappa y+2q_{02}L+4q_{04}(2LH-A).\end{array}$$
Below, we are going to express the one-form $\omega_3=s_2\omega$ in the form
$\omega_3=dS_3+s_3dH$ in order to calculate
$M_4(h)=\oint_{\delta(h)}\omega_4$ where $\omega_4=s_3\omega$. As above, we
can perform our calculations modulo exact forms. Thus,
$$
\omega_3=s_2\omega=s_2LdH+ \mbox{\rm (odd part)}\;+\; \mbox{\rm (even part)},
$$
$$
\begin{array}{rl}
\mbox{\rm (odd part)}=&
\kappa y[(2q_{20}-L)dB+d((q_{02}-\frac13L)y^2+q_{04}y^4)]\\
&+[\frac23L^2+2q_{02}L+4q_{04}(2LH-A)]dQ_1\\
&=\kappa y[(2q_{20}-L)dB-\frac13d(Ly^2)]\\
&-Q_1(\frac43LL'+2q_{02}L'+8q_{04}HL'-8q_{04}BL')dx-8q_{04}LQ_1dH\\
&=\kappa y(2q_{20}+\frac13L+2q_{02}+4q_{04}y^2)dB-\frac13\kappa yd(Ly^2)
-8q_{04}LQ_1dH\\
&=[\kappa y(2q_{20}+\frac13L+2q_{02}+4q_{04}y^2)-8q_{04}LQ_1]dH\\
&-\frac13\kappa yd(Ly^2)-\frac13\kappa y^2Ldy\\
&=[\kappa y(2q_{20}+\frac13L+2q_{02}+4q_{04}y^2)-8q_{04}LQ_1]dH.
\end{array}$$
We used  that $-Q_1L'=\kappa yB'$ and $\frac12y^2=H-B$. Similarly, by using the
identity $(2q_{20}-L)dB=d[(2q_{20}-L)B+\frac12A]$ one obtains
$$ \begin{array}{rl}
\mbox{\rm (even part)}=&\kappa ydQ_1+[\frac23L^2+2q_{02}L+4q_{04}(2LH-A)]\times\\
&\times[(2q_{20}-L)dB+d((q_{02}-\frac13L)y^2+q_{04}y^4)]\\
&=-\kappa Q_1dy-[(2q_{20}-L)B+\frac12A+(q_{02}-\frac13L)y^2+q_{04}y^4]\times\\
&\times d[\frac23L^2+2q_{02}L+4q_{04}(2LH-A)]\\
&=-[(2q_{20}-L)B+\frac12A+(q_{02}-\frac13L)y^2+q_{04}y^4]
[\frac43L+2q_{02}+4q_{04}y^2]dL\\
&-8q_{04}L[(2q_{20}-L)B+\frac12A+(q_{02}-\frac13L)y^2+q_{04}y^4]dH
-\kappa y^{-1}Q_1dH\\
&=-[2q_{02}^2+\frac23q_{02}L-\frac49L^2+4q_{04}((2q_{20}-L)B+\frac12A)]y^2dL\\
&-(6q_{02}q_{04}y^4+4q_{04}^2y^6)dL\\
&-8q_{04}L[(2q_{20}-L)B+\frac12A+(q_{02}-\frac13L)y^2+q_{04}y^4]dH
-\kappa y^{-1}Q_1dH\\
&=[4q_{02}^2L+\frac23q_{02}L^2-\frac{8}{27}L^3+8q_{04}X]dH\\
&+[24q_{02}q_{04}(2LH-A)+96q_{04}^2(LH^2-AH+Y)]dH\\
&-8q_{04}L[(2q_{20}-L)B+\frac12A+(q_{02}-\frac13L)y^2+q_{04}y^4]dH
-\kappa y^{-1}Q_1dH
\end{array}$$
where $dX=[(2q_{20}-L)B+\frac12A]dL$ and $dY=B^2dL$.
Finally, summing up all terms with $dH$, we obtain the expression
$$\begin{array}{rl}
s_3=&\kappa y(2q_{20}+\frac43L+2q_{02}+4q_{04}y^2)-8q_{04}LQ_1\\
&+4q_{02}^2L+\frac83q_{02}L^2+\frac{10}{27}L^3+4q_{04}(L+6q_{02})(2LH-A)\\
&+8q_{04}X+96q_{04}^2(LH^2-AH+Y)\\
&-8q_{04}L[(2q_{20}-L)B+\frac12A+(q_{02}-\frac13L)y^2+q_{04}y^4]
-\kappa y^{-1}Q_1.
\end{array}$$
It should be mentioned that some terms in $s_3$ were included in $s_2$ and
$s_1=L$, too. Since $M_2(h)=\oint_{\delta(h)}s_1\omega\equiv 0$,
the terms $H^kL$ have no impact in the values of
$M_3(h)=\oint_{\delta(h)}s_2\omega$ and $M_4(h)=\oint_{\delta(h)}s_3\omega$.
In the proof of Proposition \ref{p3}, we have established that
$J_2=\oint_{\delta(h)}(2LH-A)\omega\equiv 0$. By $M_3(h)\equiv 0$, one obtains
that the terms $H^kA$ and $\frac23L^2+\kappa y$ will have no impact on the
value of $M_4(h)$, too. Using these facts, one can rewrite $M_4(h)$ in the form
$$M_4(h)=\oint_{\delta(h)}(\sigma_1\omega+\sigma_2\omega+\sigma_3\omega)=
K_1+K_2+K_3$$
where
$$\begin{array}{rl}
\sigma_1=&\kappa y(2q_{20}-2q_{02}+\frac43L+4q_{04}y^2)-8q_{04}LQ_1\\
 &  -8q_{04}L[(2q_{20}-L)B+\frac12A+(q_{02}-\frac13L)y^2+q_{04}y^4],\\
\sigma_2=&\frac{10}{27}L^3+4q_{04}L(2LH-A)+8q_{04}X+96q_{04}^2Y,\\
\sigma_3=&-\kappa y^{-1}Q_1.\end{array}$$
Below, we are going to verify that $K_1+K_2=0$. Therefore
$$M_4(h)=\oint_{\delta(h)}\sigma_3\omega=\oint_{\delta(h)}\sigma_3dQ_1=
\kappa \oint_{\delta(h)} y^{-1}Q_{1x}Q_1dx.$$
Then the three formulas in Proposition \ref{p4} follow by simple calculations
making use of (\ref{I345}). Since it is assumed that
$(|\lambda|+|\mu|)q_{11}\neq 0$, $M_4(h)$ is not identically zero.
Note that the coefficient at the third formula in Proposition \ref{p4} vanishes
for $2\lambda+3\mu=0$, however this is equivalent to $a=\frac89$, a value
corresponding to the symmetric eight loop, which was excluded from
consideration here.

To finish the proof, it remains to calculate $K_2$ and $K_1$. We obtain
(modulo one-forms $dR+rdH$ which yield zero integrals)
$$\begin{array}{rl}
\sigma_2\omega&=\sigma_2dQ_1=
-Q_1d[\frac{10}{27}L^3+4q_{04}(2L^2H-LA)+8q_{04}X+96q_{04}^2Y]\\
&=-Q_1[\frac{10}{9}L^2L'+8q_{04}(2LH-\frac12A-BL)L'\\
&+8q_{04}((2q_{20}-L)B+\frac12A)L'+96q_{04}^2B^2L']dx\\
&=-Q_1L'[\frac{10}{9}L^2+8q_{04}(2q_{20}B+Ly^2)+96q_{04}^2B^2]dx\\
&=\kappa y[\frac{10}{9}L^2+8q_{04}(2q_{20}B+Ly^2)+96q_{04}^2B^2]dB\\
&=\kappa y(\frac{10}{9}L^2+8q_{04}Ly^2)d(H-\frac12y^2)
=-\kappa(\frac{10}{9}L^2y^2+8q_{04}Ly^4)dy\\
&=\kappa(\frac{20}{27}Ly^3+\frac85q_{04}y^5)L'dx.
\end{array}$$
Finally (again modulo one-forms $dR+rdH$),
$$\begin{array}{rl}
\sigma_1\omega&=
[\kappa y(2q_{20}-2q_{02}+\frac43L+4q_{04}y^2)-8q_{04}LQ_1]\omega\\
 &-8q_{04}L[(2q_{20}-L)B+\frac12A+(q_{02}-\frac13L)y^2+q_{04}y^4]dQ_1\\
 &= [\kappa y(2q_{20}-2q_{02}+\frac43L+4q_{04}y^2)-8q_{04}LQ_1]\omega\\
 &+8q_{04}Q_1L'[(2q_{20}-L)B+\frac12A+(q_{02}-\frac13L)y^2+q_{04}y^4]dx
  +8q_{04}LQ_1\omega\\
 &=\kappa y(2q_{20}-2q_{02}+\frac43L+4q_{04}y^2)\omega\\
 &-8q_{04}\kappa y[(2q_{20}-L)B+\frac12A+(q_{02}-\frac13L)y^2+q_{04}y^4]dB\\
 &=\kappa y(2q_{20}-2q_{02}+\frac43L+4q_{04}y^2)\omega\\
 &+8q_{04}\kappa y^2[(2q_{20}-L)B+\frac12A+(q_{02}-\frac13L)y^2+q_{04}y^4]dy\\
 &=\kappa y(2q_{20}-2q_{02}+\frac43L+\frac43q_{04}y^2)\omega.
\end{array}$$
Since
$$\begin{array}{rl}
\omega&=d[(2q_{20}-L)B+\frac12A+(q_{02}-\frac13L)y^2+q_{04}y^4]\\
&=[(2q_{02}-2q_{20}+\frac13L)y+4q_{04}y^3]dy-\frac13y^2dL,
\end{array}$$
we last obtain by easy calculations that
$\sigma_1\omega=-\kappa(\frac{20}{27}Ly^3+\frac85q_{04}y^5)dL$.
Proposition \ref{p4} is proved.   $\Box$

\section{The Petrov module}
We will use the notation introduced in the previous sections.

The set of Abelian integrals
$$
{\cal A}_H = \left\{ \int_{\delta(h)} \omega : \omega =  Pdx+Qdy, P, Q \in \R[x,y]\right\}
$$
is a real vector space, but also a $\R[h]$ module generated by $I_0, I_1, I_2$ with the multiplication
$$
h\cdot \int_{\delta(h)} \omega = \int_{\delta(h)} H(x,y)\omega
$$
Consider also the Petrov module
$$
{\cal P}_H = \frac{\Omega^1 }{ d \Omega^0 + \Omega^0 d H }
$$
where $\Omega^1$ is the vector space of polynomial one-forms on $\R^2$, 
and $ \Omega^0= \R[x,y]$. It is a $\R[h]$ module with multiplication
$$
h\cdot \omega = H(x,y) \omega .
$$
Let $h$ be a non-critical value of $H$. The complex algebraic curve $$\Gamma_h = \{x,y)\in \C^2: H(x,y)= h \}$$ 
has the topological type of a torus with two punctures. It follows that its first homology (co-homology) group is of
 dimension three.
According to \cite[Theorem 1.1]{gavr98} the Petrov module ${\cal P}_H$ associated  to $H$ is a free $\R[h]$ module
 generated by $\omega_0,\omega_1,\omega_2$.
\begin{proposition}[\cite{gavr98}, Proposition 3.2]
The natural map
\begin{eqnarray*}
{\cal P}_H & \rightarrow & {\cal A}_H \\
\omega & \mapsto & \int_{\delta(h)} \omega
\end{eqnarray*}
is an isomorphism of $\R[h]$ modules.
\end{proposition}

\vspace{1ex}
\noindent
{\bf Proof.}
The method of proof of the above Proposition goes back to Ilyashenko \cite{ilya69}.
The claim  follows from \cite[Proposition 3.2]{gavr98} except in the cuspidal case ($a=1$), and for the exterior 
annulus in the eight loop case.
In the case of a cuspidal loop   we note that, by making use of the Picard-Lefschetz formula, 
  the orbit of $\delta(h)$ under the monodromy group of $H$ spans the first homology group of $\Gamma_h$. 
Therefore the arguments  \cite[Proposition 3.2]{gavr98} apply.
In the case when $\delta(h)$ is represented by  an oval which belongs to the exterior period 
annulus of $\{dH=0\}$ (in the so called eight loop case), the cycle $\delta(h)$ turns out to be vanishing along 
a suitable path in the complex $h$-plane, although
 this is less obvious - see \cite[page 1170 and Fig.4]{gail05}  for a proof. The arguments used in the proof 
of \cite[Proposition 3.2]{gavr98} apply once again. $\Box$

The above result implies the following
\begin{corollary}\label{star}
Let $\alpha(h), \beta(h), \gamma(h))$ be  (real or complex) polynomials in $h$. The Abelian integral
\begin{equation}
\label{abelian}
I(h)= \alpha(h)I_0(h)+ \beta(h)I_1(h)+ \gamma(h)I_2(h)
\end{equation}
is identically zero, if and only if  $\alpha(h), \beta(h), \gamma(h)$ are identically zero.
\end{corollary}
From now on we denote by ${\cal A}_n$ the space of of Abelian integrals of the form (\ref{abelian}), with
$$\deg \alpha \leq n, \deg \beta \leq n, \deg \gamma \leq n.$$
\begin{corollary}
The maximal dimension of the vector space ${\cal A}_n$
equals $3(n+1)$.
\end{corollary}
\begin{remark}
\label{remark1}
The vector space of Abelian integrals ${\cal A}_n$   coincides with the space of Abelian integrals
\begin{equation}
\label{abelianp}
\int_{\delta(h)} P(x,y) dx + Q(x,y) dy
\end{equation}
where $P,Q$ are real polynomials of weighted degree $4n+5$, where the weight of $x$ is $1$ and the weight 
of $y$ is $2$. Therefore the result follows also from  {\em \cite[p.582]{gavr98}}.
\end{remark}

\section{Zeroes of Abelian integrals}
In this section we find upper bounds for the number of the zeroes of the Abelian integrals  ${\cal A}_n$ defined
 in (\ref{abelian}) on the interval of existence of the ovals $\delta(h)$.
Similar results were earlier obtained for the space of non-weighted Abelian integrals (\ref{abelianp})
 ($\deg P, \deg Q \leq n$) by Petrov \cite{petr90} and Liu \cite{liu03}, see the survey of Christopher and Li \cite{chli07}.
 Our upper bounds however do not follow from the aforementioned papers, see Remark \ref{remark1}. 
They  will be proved  along the  lines, given in \cite[section 7]{gavr99}.

All families of cycles will depend continuously on a parameter $h$ and will be defined without ambiguity
 in the complex half-plane $\{h: Im (h) >0\}$.
This will allow a continuation on $\C$ along any curve avoiding the real critical values of $H$.
 In particular, it will be supposed that all three critical values of $H$ are real.

\begin{figure}
\begin{center}
 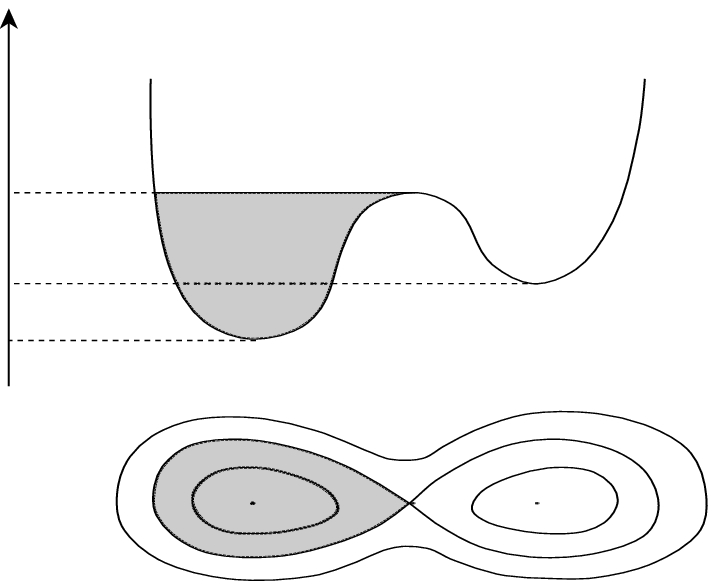
 \end{center}
 \caption{The graph of the polynomial $\frac12 x^2-\frac23 x^3 +\frac{a}{4} x^4 $,
 $\frac89<a<1$,
 and the level sets $\{H=h\}$}
 \label{fig1}
 \end{figure}

\subsection{The interior eight-loop case}
Using the normal form (\ref{ham}) we can suppose that $\frac89 < a < 1$.
Let $\delta (h) \subset \{H=h\}$ be a
continuous family of ovals defined on a maximal open interval
$\Sigma =(h_c,h_s)$, where for $h=h_c=0$ the oval degenerates to a point
$\delta (h_c)$ which is a center and for $h=h_s>0$ the oval becomes a
homoclinic loop of the Hamiltonian system $dH=0$. The family
$\{ \delta (h)\}$ represents a continuous family of cycles vanishing
at the center $\delta (h_c)$.
\begin{theorem}
The space of Abelian integrals  ${\cal A}_n$ corresponding to Fig. \ref{fig1} is Chebyshev on the interval of existence of $\delta(h)$.
\end{theorem}

\begin{figure}
 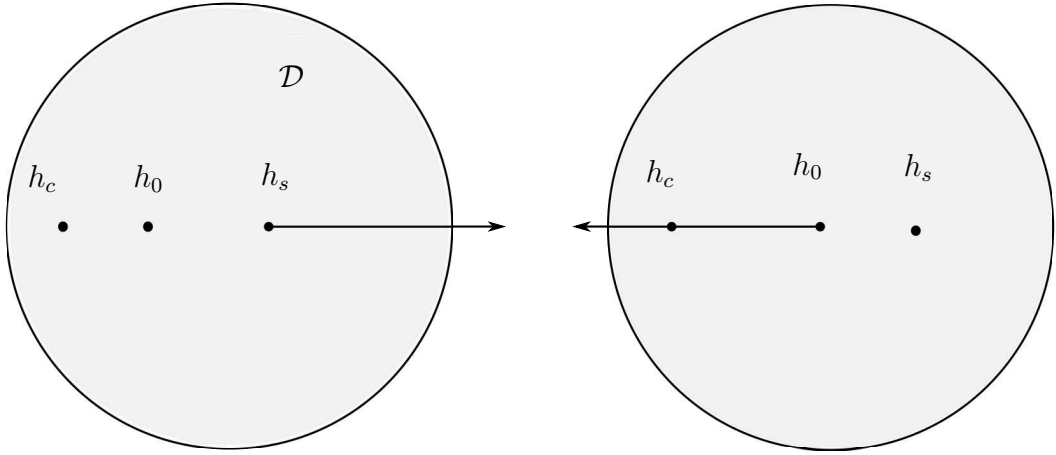
 \caption{The branch cuts and the domains of the Abelian integral $I(h)$ and $ W_{\delta,\delta_s}(\omega',\omega_0') $ respectively.}
 \label{fig1plus}
 \end{figure}

\vspace{1ex}
\noindent
{\bf Proof.}
We note first  that $I_0(h), I_1(h), I_2(h)$  can be expressed as linear combinations of
 $I_0'(h), I_1'(h), I_2'(h)$, whose coefficients are polynomials in $h$ of degree one. Therefore the vector space
$$
{\cal A}_n' = \{I'(h): I(h) \in {\cal A}_n \}
$$
coincides with the vector space of Abelian integrals
$$
\{ \alpha(h)I_0'(h)+ \beta(h)I_1'(h)+ \gamma(h)I_2'(h) : \deg \alpha \leq n, \deg \beta \leq n, \deg \gamma \leq n \} .
$$
We shall prove the Chebyshev property of ${\cal A}_n' $ in the complex domain
$$
{\cal D } = \C \setminus [h_s, \infty) .
$$
in which $I'(h)$ has an analytic extension, see fig.\ref{fig1plus}. For this purpose we apply the argument principle to the function
$$
F(h) = \frac{I'(h)}{I_0'(h)}.
$$

We note that $I_0'(h)$ is a complete elliptic integral of first kind and hence cannot vanish in ${\cal D }$. For sufficiently big $|h|$
 the function $F(h)$ behaves as $h^{n+\frac12}$ and hence the increment of the argument of $F$ along a circle with a sufficiently
 big radius is close to $(2n+1)\pi$.
Along the interval $[h_s, \infty)$ the imaginary part of $F(h)$ can be computed by making use of the Picard-Lefshetz formula.
 Namely, let $\{\delta_s(h)\}_h$ be a continuous family of cycles, vanishing at the saddle point as $h$ tends to $h_s$. Then along
 $[h_s, \infty)$ the family $\delta(h)$ has two analytic complex-conjugate continuations $\delta^\pm (h)$, $\delta^+ = \bar{\delta^-}$ 
and moreover, by the Picard-Lefshetz formula the cycle
$$\delta^+(h) - \delta^-(h)= \delta_s(h)
$$
where the identity should be understood up to homology equivalence. This implies the following identity along
$[h_s, \infty)$
$$
2 Im(F(h)) = \frac{\int_{\delta^+(h)} \omega'}{\int_{\delta^+_0(h)} \omega_0'}-
 \frac{\int_{\delta^-(h)} \omega'}{\int_{\delta^-_0(h)} \omega_0'}
 =
 \frac{W_{\delta,\delta_s}(\omega',\omega_0')}{|I'_0(h)|^2}
$$
where
$$
W_{\delta,\delta_s}(\omega',\omega_0') =
\det
\left(
\begin{array}{ccc}
 \int_{\delta(h)} \omega' &  \int_{\delta_s(h)} \omega'   \\
 \int_{\delta(h)} \omega'_0 &  \int_{\delta_s(h)} \omega'_0
\end{array}
\right) .
$$
Following \cite[section 7]{gavr99} we may use the reciprocity law on the elliptic curve $\{H=h\} $ to compute
$$
W_{\delta,\delta_s}(\omega',\omega_0') =
p(h)+ q(h)\int_{-\infty}^{+\infty} \frac{dx}{y}
$$
where $p(h),q(h)$ are suitable degree $n$ polynomials, $ \pm \infty$ are the two "infinite" points on the compactified Riemann surface
 $\Gamma_h$, and the integration is along some path connecting  $ \pm \infty$ on $\Gamma_h$.

It is easy to check now that the function $p(h)+ q(h)\int_{-\infty}^{+\infty} \frac{dx}{y}$ can have at most $2n+1$ zeroes on
 $[h_s, \infty)$. For this, consider an analytic continuation of this function to the complex domain $\C \setminus (-\infty,h_0]$
 where $h_0$ is a critical value of $H$, $h_0< h_s$, see fig.\ref{fig1plus}. By the Picard-Lefshetz formula, the imaginary part
 of $p(h)+ q(h)\int_{-\infty}^{+\infty} \frac{dx}{y}$ along the branch cut $(-\infty,h_0)$ equals
$$
\tilde{q}(h)\int_{\delta_(h)}  \frac{dx}{y}
$$
where $\tilde{q}$ differs from $q$ by an addition of a constant.
 We conclude that the imaginary part of this function vanishes at most $n$ times. This combined to the asymptotic behavior
$$
p(h)+ q(h)\int_{-\infty}^{+\infty} \frac{dx}{y} \sim h^{n}\times const.
$$
gives that the increase of the argument along a big circle is close to $2\pi n$ and finally, that our function can have at most
 $2n+1$ zeroes on $\C \setminus (-\infty,h_0]$. Now we come back to the function $F(h)$ and conclude that it can have at
 most $3n+2$ zeroes in the complex domain ${\cal D}$, counted with the multiplicity. As $I(0)=0$ the same conclusion holds 
true for $I(h)$ on the real interval $(-\infty, h_s)$. $\Box$

\begin{remark} Through the proof we did not inspect the behavior of $F(h)$ near the branch point $h_s$. In the original papers
 of Petrov a small circle centered at $h_s$ is removed and the behavior of $F$ along it is taken into account. It is important to note
 that, we do not remove a small circle here, because we use a slightly improved version of the argument principle,  as explained in 
section 2.4 of \cite{gavr13}. It allows one to apply the argument principle, even if $F$ is not analytic at $F(h_s)$, provided that
 $F$ has a continuous limit at $h_s$, which is not zero. The case when $F(h_s) = 0$ is studied then by a small perturbation 
(by adding a real constant for instance) - this does not decrease the number of zeroes of $F$ in the complement of the branch cut.
Of course, the same considerations hold true for the function $\int_{-\infty}^{+\infty} \frac{dx}{y}$ in its respective domain 
of analyticity.\end{remark}

\subsection{The saddle-loop case}

\begin{figure}
\begin{center}
 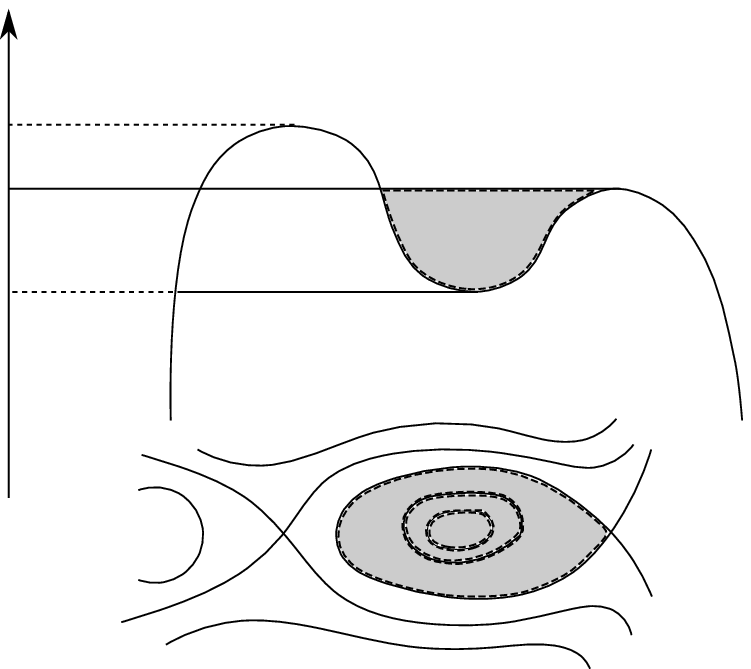
 \end{center}
 \caption{The graph of the polynomial $\frac12 x^2-\frac23 x^3 +\frac{a}{4} x^4 $,
 $a<0$,
 and the level sets $\{H=h\}$}
 \label{fig3}
 \end{figure}

In the normal form (\ref{ham}) we suppose that $a < 0$.
As before, we let $\delta (h) \subset \{H=h\}$ be a
continuous family of ovals defined on a maximal open interval
$\Sigma =(h_c,h_s)$, where for $h=h_c=0$ the oval degenerates to a point
$\delta (h_c)$ which is a center and for $h=h_s>0$ the oval becomes a
homoclinic loop of the Hamiltonian system $dH=0$. The family
$\{ \delta (h)\}$ represents a continuous family of cycles vanishing
at the center $\delta (h_c)$.
\begin{theorem}
The space of Abelian integrals  ${\cal A}_n$ corresponding to the shadowed area on Fig. \ref{fig3} is of dimension $3n+3$, 
and each Abelian integral from ${\cal A}_n$ can have at most $4n+3$ zeroes.
\end{theorem}

\vspace{2ex}
\noindent
{\bf Proof.}
We shall prove the Chebyshev property of ${\cal A}_n' $ in the complex domain
$$
{\cal D } = \C \setminus [h_s, \infty) .
$$
in which $I'(h)$ has an analytic extension.  For this purpose we apply the argument principle to the function
$$
F(h) = \frac{I'(h)}{I_0'(h)}.
$$
Indeed, a local analysis shows that at $h_s, h_0$ the function $F|_{\cal D } $ has   continuous limits, which we assume 
to be non-zero. $I_0'(h)$ is a complete elliptic integral of first kind and hence cannot vanish in ${\cal D }$. For sufficiently 
big $|h|$ the function $F(h)$ behaves as $h^{n+\frac12}$ and hence the increment of the argument of $F$ along a circle
 with a sufficiently big radius is close to $(2n+1)\pi$.
Along the intervals $(h_s, h_0)$ and $(h_0, \infty)$ the imaginary part of $F(h)$ can be computed by making use of the
 Picard-Lefshetz formula. Namely, let $\{\delta_s(h)\}_h$, $\{\delta_0(h)\}_h$ be the continuous family of cycles, 
vanishing at the saddle points $h_s$ and $h_0$ respectively, as $h$ tends to $h_s$
and $h_0$. As in the preceding section we deduce that along $[h_s, h_0)$,
$$
2 Im(F(h))
 =
 \frac{W_{\delta,\delta_s}(\omega',\omega_0')}{|I'_0(h)|^2}, \;\; h \in (h_s, h_0)
$$
while along $ (h_0, \infty)$
$$
2 Im(F(h))  =
 \frac{W_{\delta,\delta_s}(\omega',\omega_0')}{|I'_0(h)|^2}
 +
 \frac{W_{\delta,\delta_0}(\omega',\omega_0')}{|I'_0(h)|^2}
 , \;\; h \in (h_0, \infty) .
$$
The function
$$W_{\delta,\delta_s}(\omega',\omega_0'), \; h \in (h_s, h_0)$$
  allows an analytic continuation in $\C \setminus [h_0,\infty)$ and exactly as in the preceding section we compute that it 
can have at most $2n+1$ zeroes there. More precisely, $W_{\delta,\delta_s}(\omega',\omega_0')$
has an  analytic continuation in $\C \setminus [h_0,\infty)$. The number of its zeroes in this domain is bounded by $n$ 
(coming from the behavior at infinity) plus one plus the number of the zeroes of
$$
2 Im ( W_{\delta,\delta_s}(\omega',\omega_0') )= W_{\delta_0,\delta_s}(\omega',\omega_0')
 = q(h)\int_{-\infty}^{+\infty} \frac{dx}{y}, \;\; h\in (h_0,\infty)
$$
where $q$ is a degree $n$ polynomial.
Similarly, the function
$$W_{\delta,\delta_s}(\omega',\omega_0') + W_{\delta,\delta_0}(\omega',\omega_0'), \;  (h_0, \infty)$$
allows an analytic continuation in $\C \setminus [h_s,h_0]$ and its zeroes there are bounded by $n$ plus plus one plus
 the number of the zeroes of
$$
2 Im( W_{\delta,\delta_s}(\omega',\omega_0') + W_{\delta,\delta_0}(\omega',\omega_0'))=
W_{\delta_s,\delta_0}(\omega',\omega_0')= q(h)\int_{-\infty}^{+\infty} \frac{dx}{y}, \;\; h\in (h_s,h_0) .
$$
Summing up the above information, we get that
 the function $F(h)$ can have at most $4n+3$
 zeroes in the complex domain ${\cal D}$, counted with the multiplicity. As $I(0)=0$ the same conclusion holds true
 for $I(h)$ on the real interval $(-\infty, h_s)$. $\Box$

\subsection{ The exterior eight-loop case}

\begin{figure}
\begin{center}
 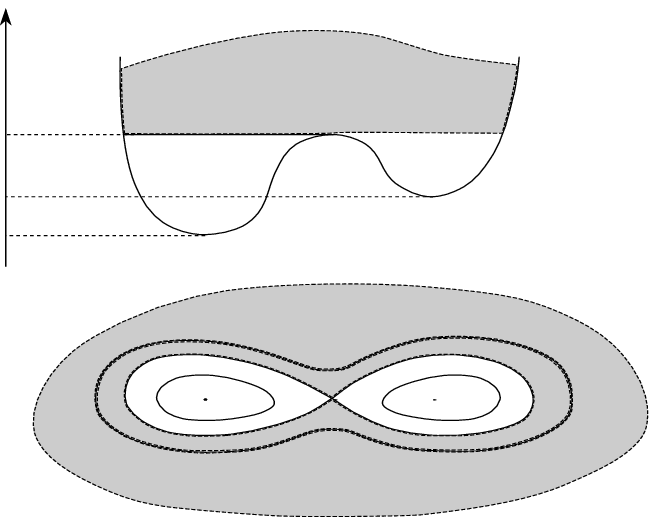
 \end{center}
 \caption{The graph of the polynomial $\frac12 x^2-\frac23 x^3 +\frac{a}{4} x^4 $,
 $\frac89<a<1$,
 and the level sets $\{H=h\}$}
 \label{fig4}
 \end{figure}

In this section we consider the exterior eight-loop case, with period annulus as shown on fig.\ref{fig4} and
$\frac89 < a < 1$.
Let $\delta (h) \subset \{H=h\}$ be the
continuous family of ovals defined on the maximal open interval
$\Sigma =(h_s,\infty)$
\begin{theorem}
The space of Abelian integrals  ${\cal A}_n$ corresponding to the shadowed area on Fig. \ref{fig4} is of dimension $3n+3$, 
and each Abelian integral from ${\cal A}_n$ can have at most $4n+4$ zeroes.
\end{theorem}

\vspace{2ex}
\noindent
{\bf Proof.}
We shall evaluate the number of the zeroes of a function from ${\cal A}_n' $ in the complex domain
$$
{\cal D } = \C \setminus (-\infty,h_s] .
$$
in which $I'(h)$ has an analytic extension.  For this purpose we apply the argument principle to the function
$$
F(h) = \frac{I'(h)}{I_0'(h)}.
$$
As before, a local analysis shows that at $h_s, h_c^1, h_c^2$ the function $F|_{\cal D } $ has   continuous limits, which 
we assume to be non-zero. $I_0'(h)$ is a complete elliptic integral of first kind and hence 
cannot vanish in ${\cal D }$. For sufficiently big $|h|$ the function $F(h)$ behaves as $h^{n+\frac12}$ and hence the 
increment of the argument of $F$ along a circle with a sufficiently big radius is close to $(2n+1)\pi$.
It remains to study the number of the zeroes of the imaginary part of $F(h)$ along the intervals
$$
(-\infty , h_c^1), \; (h_c^1,h_c^2), \; (h_c^2,h_s) .
$$
Namely, let $\{\delta_s(h)\}_h$, $\{\delta_c^1(h)\}_h$, $\{\delta_c^2(h)\}_h$ , where $Im h \geq 0$, be the continuous 
family of cycles, vanishing at the saddle points as $h$  tends to $h_s$, and $h_c^1$ or $h_c^2$, respectively. These cycles
 are defined
up to an orientation, and we consider their continuation to ${\cal D } = \C \setminus (-\infty,h_s]$, as well the limits
 along the branch cut $(-\infty,h_s]$. The family of exterior loops $\{\delta(h)\}$ is expressed in terms of these vanishing cycles as
follows
$$
\delta(h) = \delta_c^1(h) + \delta_c^2(h) + \delta_s(h), \; h \in {\cal D }
$$
(the orientations of the vanishing cycles are fixed from this identity). Let $\delta^+(h)=\delta(h)$ be the continuation of $\delta(h)$
 on $(-\infty,h_s]$, along paths contained in the upper complex half-plane, and $\delta^-(h)$ be the continuation on $(-\infty,h_s]$ 
along paths contained in the lower complex half-plane. The Picard-Lefschetz formula easily implies
$$
\delta^-(h) = \delta_c^1(h) + \delta_c^2(h) - \delta_s(h), h \in (h_c^2,h_s)
$$
$$
\delta^-(h) = \delta_c^1(h)  - \delta_s(h), h \in (h_c^1,h_c^2)
$$
$$
\delta^-(h) = - \delta_s(h), h \in (-\infty, h_c^1)
$$

As in the preceding section we deduce that along the branch cut $ (-\infty,h_s)$
we have
\begin{equation}
\label{im1}
2 Im(F(h))
 =
 \frac{W_{\delta,2\delta_s}(\omega',\omega_0')}{|I'_0(h)|^2}, \;\; h \in (h_c^2, h_s)
\end{equation}
and
\begin{equation}
\label{im2}
2 Im(F(h))  =
 \frac{W_{\delta,2\delta_s}(\omega',\omega_0')}{|I'_0(h)|^2}
 +
 \frac{W_{\delta,\delta_c^2}(\omega',\omega_0')}{|I'_0(h)|^2}
 , \;\; h \in (h_c^1, h_c^2)
\end{equation}
and
\begin{equation}
\label{im3}
2 Im(F(h))  =
 \frac{W_{\delta,2\delta_s}(\omega',\omega_0')}{|I'_0(h)|^2}
 +
 \frac{W_{\delta,\delta_c^1+\delta_c^2}(\omega',\omega_0')}{|I'_0(h)|^2}= 
\frac{W_{\delta,\delta_s}(\omega',\omega_0')}{|I'_0(h)|^2}, \;\; h \in (-\infty,h_c^1) .
\end{equation}
Clearly, the function $W_{\delta,\delta_s}(\omega',\omega_0')$ has an analytic continuation in 
$\C \setminus [h_c^1, h_c^2]$. Its number of zeroes in this domain depends on the zeroes of
$$
2 Im( W_{\delta,\delta_s}(\omega',\omega_0')) = W_{\delta_c^1,\delta_c^2}(\omega',\omega_0')
=q(h)\int_{-\infty}^{+\infty} \frac{dx}{y}, \;\; h\in (h_c^1, h_c^2) .
$$
Thus, the total number of the zeroes of the functions (\ref{im1}), (\ref{im3}) is bounded by $n+1$ plus the 
number of the zeroes of $q(h)$ on the interval
$(h_c^1, h_c^2)$. Finally, similar considerations show that the function (\ref{im2}) has an analytic continuation in
$$
\C \setminus \{ (-\infty, h_c^1) \cup (h_c^2, \infty) \}
$$
and its zeroes in this domain are bounded by $n+1$ plus the number of the zeroes of the polynomial $q(h)$ 
on the interval $ (-\infty, h_c^1) \cup (h_c^2, \infty) $.

Summing up the above information, we get that
 the function $F(h)$ can have at most $4n+3$
 zeroes in the complex domain ${\cal D}$, counted with the multiplicity. Therefore the Abelian integral $I(h)$ 
has at most $4n+4$ zeroes on the real interval $(-\infty, h_s)$.  $\Box$

\section{Lower bounds for the number of zeroes of $M_k(h)$}

In this section we provide examples which show that Chebyshev's property would not hold
in the saddle-loop case. For this purpose, we study the number of small-amplitude limit cycles 
bifurcating around the center at the origin. 

We begin with the system satisfied by the basic integrals $I_k(h)$. It is
derived in a standard way by using (\ref{ham}), (\ref{I345}) and the formula
$I_k'(h)=\oint_{\delta(h)}(x^k/y)dx$.

\begin{lemma}\label{l1}The integrals $I_0(h)$, $I_1(h)$ and $I_2(h)$
satisfy the system
$$\begin{array}{l}
\frac43hI_0'-\frac{2}{9a}I_1'-(\frac13-\frac{4}{9a})I_2'=I_0,\\
\frac{2}{9a}hI_0'+(h+\frac{1}{4a}-\frac{10}{27a^2})I_1'
-(\frac{13}{18a}-\frac{20}{27a^2})I_2'=I_1,\\
-(\frac{4}{15a}-\frac{56}{135a^2})hI_0'+(\frac{4}{15a}h+\frac{29}{45a^2}-\frac{56}{81a^3})I_1'
+(\frac45h+\frac{4}{15a}-\frac{46}{27a^2}+\frac{112}{81a^3})I_2'=I_2.
\end{array}$$
\end{lemma}

\noindent
We use this system to find the expansions of integrals $I_k$, $k=0,1,2$
near $h=0$. Denoting $c=I'_0(0)\neq 0$, one obtains

\begin{lemma}\label{l2}
The following expansions hold near $h=0$:
$$\begin{array}{rl}
I_0(h) & =c[h+(\frac53-\frac38a)h^2+(\frac{385}{27}-\frac{35}{4}a
+\frac{35}{64}a^2)h^3\\[2mm]
& + (\frac{85085}{486}-\frac{25025}{144}a+\frac{5005}{128}a^2
-\frac{1155}{1024}a^3)h^4\\[2mm]
&+1001(\frac{7429}{2916}-\frac{2261}{648}a
+\frac{1615}{1152}a^2-\frac{85}{512}a^3+\frac{45}{16384}a^4)h^5+\ldots],\\[2mm]
I_1(h) & =c[h^2+(\frac{70}{9}-\frac{35}{12}a)h^3
+(\frac{5005}{54}-\frac{5005}{72}a+\frac{1155}{128}a^2)h^4\\[2mm]
&+1001(\frac{323}{243}-\frac{323}{216}a+\frac{85}{192}a^2-\frac{15}{512}a^3)h^5
\\[2mm]
&+1001(\frac{185725}{8748}-\frac{185725}{5832}a+\frac{52003}{3456}a^2-
\frac{11305}{4608}a^3+\frac{1615}{16384}a^4)h^6 +\ldots],\\[2mm]
I_2(h) & =c[\frac12h^2+(\frac{35}{9}-\frac58a)h^3
+(\frac{5005}{108}-\frac{385}{16}a+\frac{315}{256}a^2)h^4\\[2mm]
&+1001(\frac{323}{486}-\frac{85}{144}a+\frac{15}{128}a^2-\frac{3}{1024}a^3)h^5
\\[2mm]
&+1001(\frac{185725}{17496}-\frac{52003}{3888}a+\frac{11305}{2304}a^2
-\frac{1615}{3072}a^3+\frac{255}{32768}a^4)h^6 +\ldots].
\end{array}$$
\end{lemma}

\vspace{2ex}
\noindent
{\bf Proof.} We rewrite the system from Lemma \ref{l1} in the
form $({\bf A}h+{\bf B}){\bf I}'={\bf E}{\bf I}$ where
${\bf I}=(I_0,I_1,I_2)^\top$. As ${\bf I}(h)$ is a solution which is
analytical near zero and ${\bf I}(0)=0$, one can replace
$${\bf I}(h)=\sum_{k=1}^\infty {\bf V}_k h^k,\quad
{\bf V}_k=(V_{0,k},V_{1,k},V_{2,k})^\top$$
in the system. Then the coefficient at $h^k$ should be zero, which yields
the equation
$$(k+1){\bf B}{\bf V}_{k+1}=({\bf E}-k{\bf A}){\bf V}_k.$$
Since ${\bf I}'(0)=(c,0,0)^\top= {\bf V}_1$, one can solve the system above
with respect to $(V_{0,k},V_{1,k+1}, V_{2,k+1})$ and thus to obtain
via recursive procedure formulas for all ${\bf V}_k$, $k=2,3,\ldots$.
Explicitly,
$$\begin{array}{rl}
(8-9a)V_{1,k+1}&=[8-9a+(88-87a)\frac{k-1}{k+1}]V_{0,k}
-(48a-36a^2)\frac{k-1}{k+1}V_{1,k},\\
(8-9a)V_{2,k+1}&=[4-\frac92a+(44-\frac{63}{2}a)\frac{k-1}{k+1}] V_{0,k}
-24a\frac{k-1}{k+1}V_{1,k},\\
V_{0,k+1}&=\frac{6k-1}{3k}V_{1,k+1}-a\frac{4k-1}{4k}V_{2,k+1},\;\;\;
k=1,2,3,\ldots
\end{array}$$
Applying these formulas, we obtain the expansions in Lemma 2.   $\Box$

\vspace{2ex}
\noindent
{\bf Proof of Theorem \ref{t5}.} Consider the following linear combinations
$$\begin{array}{ll}
J_0=I_0,  & \hspace{10mm} J_3=\alpha_1 hI_0+\beta_1 I_1+\gamma_1 I_2,\\
J_1=I_1,   &\hspace{10mm} J_4=\alpha_2 hI_0+(\beta_2+\delta_2h)I_1+\gamma_2 I_2,\\
J_2=I_1-2I_2,  &\hspace{10mm} J_5=\alpha_3 hI_0+(\beta_3+\delta_3h)I_1+(\gamma_3+\eta_3 h)I_2,
\end{array}$$
where
$$\begin{array}{ll}
\alpha_1=a,   
& \hspace{10mm} \delta_2=\frac83a^2+a^3, \\[2mm]
\beta_1=-\frac{11}{3}+\frac{21}{40}a,  
  & \hspace{10mm}\alpha_3=\frac{17}{81}a-\frac{775}{5148}a^2+\frac{63}{9152}a^3, \\[2mm]
\gamma_1=\frac{22}{3}-\frac{61}{20}a, 
 & \hspace{10mm} \beta_3=-\frac{187}{243}+\frac{55}{72}a-\frac{1085}{9152}a^2-\frac{189}{73216}a^3,   \\[2mm]
\alpha_2=\frac{208}{63}a-\frac23a^2,   
 & \hspace{10mm}\gamma_3=\frac{374}{243}-\frac{631}{324}a+\frac{155}{288}a^2-\frac{315}{36608}a^3,   \\[2mm]
\beta_2=-\frac{2288}{189}+\frac{52}{9}a+\frac14a^2,  
  & \hspace{10mm}\delta_3=\frac{119}{702}a^2-\frac{147}{1144}a^3-\frac{189}{18304}a^4, \\[2mm]
\gamma_2=\frac{4576}{189}-\frac{1144}{63}a+\frac56a^2,  
 & \hspace{10mm} \eta_3=\frac{49}{234}a^3. 
\end{array}$$
The coefficients above are chosen so that $J_k(h)=O(h^{k+1})$ near zero
for $0\leq k \leq 5$. Their explicit values are determined from the respective
linear systems. By calculation, then one obtains
$$\begin{array}{ll}
J_0=c[h+\ldots],    
 & \hspace{5mm} J_3=c[\frac{49}{32}a^2(a+\frac83)h^4+(\frac{68992}{405}+O(a+\frac83))h^5+\ldots], \\[2mm]
J_1=c[h^2+\ldots], 
 & \hspace{5mm} J_4=c[\frac{154}{9}a^4h^5+\ldots],   \\[2mm]
J_2=c[-\frac53ah^3+\ldots],   
 & \hspace{5mm} J_5=c[\frac{49}{128}a^5(a+\frac89)h^6+(-119(\frac23)^{14}+O(a+\frac89))h^7+\ldots]. \\[2mm]
\end{array}$$
Let us fix the Hamiltonian parameter $a$ be a little bit smaller that $-\frac83$, so that we would have
$J_3=c[\delta_4h^4+\delta_5h^5+O(h^6)]$ with $|\delta_4|<\!\!<|\delta_5|$ and $\delta_4<0<\delta_5$.
Then, one can choose a linear combination $J(h)$ of $J_k$, $0\leq k\leq 3$, such that 
$J(h)=c\sum_{k=1}^5\delta_kh^k[1+O(h)]$ will satisfy $\delta_k\delta_{k+1}<0$ and 
$|\delta_k|<\!\!<|\delta_{k+1}|$. Therefore, $J(h)$ would have 4 small positive zeroes. 
As the four coefficients in (\ref{m1}) are independently free, 
one can take a small perturbation such that $M_1(h)=J(h)$ will produce 4 small-amplitude limit cycles around 
the center at the origin. The proof of the claim concerning $M_2(h)$ is the same, as long as we fix the
parameter $a$ a little bit smaller than $-\frac89$ and construct in the same way a linear combination
$J(h)=c\sum_{k=1}^7\delta_kh^k[1+O(h)]$ with coefficients having the same properties, thus $M_2(h)$ producing 
6 small positive zeroes in the saddle-loop case.  

For all other $a\in\R$ different from $0$, $-\frac83$ and $\pm\frac89$, any linear combination of $J_k$,
$0\leq k\leq m$ where $m=3,4,5$, will have at most $m$ small positive zeroes. Moreover, $M_k(h)$, $k=1,2,3$ can be 
expressed as linear combination of the respective $J_k$, thus having as much zeroes at its dimension minus one.
 It is easy to see that $M_4(h)$ has no small positive zeroes at all. $\Box$

\section*{Acknowledgements}
Part of the paper was written while the second author was visiting the Institut of Mathematics, University of Toulouse III (Paul Sabatier). He is obliged for the hospitality.

\def\cprime{$'$} \def\cprime{$'$} \def\cprime{$'$} \def\cprime{$'$}
  \def\cprime{$'$} \def\cprime{$'$} \def\cprime{$'$}

\end{document}